
\baselineskip=14pt
\parskip=10pt

\font\eightrm=cmr8 

\magnification=\magstephalf

\def\1{{\overline{1}}}
\def\2{{\overline{2}}}
\parindent=0pt
\overfullrule=0in

\def\frac#1#2{{#1 \over #2}}

\bf
\centerline
{
There are
${\bf \frac{1}{30}(r+1)(r+2)(2r+3)(r^2+3r+5)}$ Ways For the Four Teams 
of a World Cup Group 
}
\centerline
{
to Each Have $r$ Goals For
and $r$ Goals Against
}
\centerline
{
[Thanks to the Soccer Analog of Prop. 4.6.19 of Richard Stanley's (Classic!) EC1]
}
\rm
\bigskip
\centerline
{\it By Shalosh B. EKHAD and Doron ZEILBERGER}
\bigskip
\qquad \qquad \qquad \qquad 
{\it Dedicated to Richard Peter Stanley who just turned 
``number of ways for a simple 1D Drunkard to return home after 8 steps''-years-old}
\bigskip

{\bf Proof of the Statement in the Title}: Calling this quantity $S_4(r)$, 
and using the case $n=4$ of the Soccer analog of Prop. 4.6.19 of [EC1]
(see Comment 1 below), we see that

$\bullet$ $S_4(r)$ is a polynomial of degree $5$  \quad .

$\bullet$ $S_4(-1)=0$ \quad , \quad $S_4(-2)=0$ \quad .

$\bullet$ $S_4(-3-r) \, = \, -S_4(r)$ \quad .

Hence, it suffices to check the statement at the {\bf two} values $r=0$ and $r=1$. Obviously $S_4(0)=1$
(there is just one all-zero $4 \times 4$ matrix), and almost obviously, $S_4(1)=9$
(the number of derangements of length $4$ equals $9$). QED!

{\bf Comments} 

{\bf 1}. The same method of proof that guru Richard Stanley used to prove Prop. 4.6.19 of [EC1] 
(first edition; it is Prop. 4.6.2 in the second edition.)
[{\eightrm the first edition is missing $(-1)^{n-1}$ in front of $H_n(r)$}]
yields:

{\bf The Soccer analog of Prop. 4.6.19 of EC1}

Let $S_n(r)$ be the number of ways $n$ Soccer teams, playing a round-robin tournament, each scored a total
of $r$ ``Goals For'' (GF), and a total of $r$ ``Goals Against'' (GA), (in other words the number of $n \times n$
magic squares whose diagonal entries are all $0$ (no team plays against itself!)).

For each fixed integer $n \geq 3$, the function $S_n(r)$ is a polynomial in $r$ of degree 
$n^2-3n+1$. Since it is a polynomial, it can be evaluated at {\it any} (not necessarily positive) integer, and we have
$$
S_n(-1)=S_n(-2)= \dots = S_n(-n+2)=0 \quad ,
$$
$$
S_n(-(n-1)-r)=-S_n(r) \quad .
$$

{\bf 2.} We got the idea for this tribute to our beloved guru when we attended the Stanley@70 conference,
held at MIT, June 23-27, 2014, at the same time as the 
preliminary Group stage of the World Cup 2014.

{\bf 3.} While the proof of the statement of the title did not require any computer, the analogous result for
$n=5$ is already beyond the scope of mere humans (or the human would have to be very stupid to
spend time on it). We have
$$
S_5(r)={\frac {1}{241920}}\, \left( r+1 \right)  \left( r+2 \right)  \left( r+3 \right)  \cdot
$$
$$
\left( 43\,{r}^{8}+688\,{r}^{7}+4934\,{r}^{6}+20680\,{r}^{5}+55907\,{r}^{4
}+101272\,{r}^{3}+123436\,{r}^{2}+96240\,r+40320 \right)  \quad .
$$
Of course $S_3(r)=r+1$ (why?).

To see $S_6(r)$, please go to: {\tt http://www.math.rutgers.edu/\~{}zeilberg/tokhniot/oGOALS1} $\,\,$ .

One can get this, and (potentially!) infinitely more results, using the Maple package {\bf GOALS} available directly from:
{\tt http://www.math.rutgers.edu/\~{}zeilberg/tokhniot/GOALS} $\,\,$,
or via the front of this article:
{\tt http://www.math.rutgers.edu/\~{}zeilberg/mamarim/mamarimhtml/worldcup.html}$\,\,$.

In particular, procedure {\tt MagicPolAG(n,r,Sr,Sc,B,A)} can find polynomial expressions for the number of 
$n \times n$ matrices with non-negative integer entries, whose
row-sums are 
$$
(r+Sr[1], \dots, r+Sr[n]) \quad ,
$$ 
and whose column sums are
$$
(r+Sc[1], \dots, r+Sc[n]) \quad ,
$$ 
for {\it any} fixed numerical vectors {\tt Sr}, {\tt Sc}, of length $n$ (of course they have to add-up to the same number),
and any assignment $A$ where some entries are fixed beforehand, where $B$ denotes ``wild-card''. See the on-line help there.
See the above-mentioned front of this article  for numerous sample input and output files.

{\bf 4.} The sequence for $S_4(r)$ is already in Neil Sloane's OEIS, (see [OEIS1]), but for {\bf different} reasons!
Can you find a bijection? On the other hand, the sequence for $S_5(r)$ is not (yet!) there:
$$
 1, 44, 870, 9480, 68290, 365936, 1573374, 5709120, 18107760, 51488800, 133748186, 321979164,  \dots \quad ,
$$
but we are sure that very soon it will!

{\bf 5.} Another, even more useful, Maple package accompanying this article is  {\tt WorldCup}, 
one of whose many procedures is `{\tt Ptor}', that finds all the possible scenarios that lead to a given
{\it score board}, consisting of the ``Goals For'', ``Goals Against'', and ``PTS'', vectors.
(Recall that a win yields $3$ points, a draw, $1$ point, and a loss, $0$ points.).

The number of possible scenarios for groups A-H were as follows

{\bf 2014}: $8, 32, 7, 3, 13, 3, 12, 3$ \quad .

{\bf 2010}: $2, 3, 2, 2, 3, 3, 1, 2$ (in particular, the score-board for Group $G$ {\it uniquely} determined the individual game scores).

{\bf 2006}:  $11, 1, 5, 3, 2, 6, 1, 8$ (in particular, the score-boards for Groups $B$ and $G$ {\it uniquely} determined the individual game scores).

{\bf 2002}: $1, 12, 5, 5, 1, 2, 3, 1$ (in particular, the score-boards for Groups $A,E,H$ {\it uniquely} determined the individual game scores).

{\bf 1998}:  $9, 6, 4, 3, 3, 3, 1, 1$ (in particular, the score-boards for Groups $G,H$ {\it uniquely} determined the individual game scores).

{\bf 6.} The Israeli daily {\it Yedioth Ahronoth} published, at the start of the 2014 World Cup,  a ``soccer puzzle'' where
the above-mentioned vectors, GF, GA, and PTS, are given, and the solver has to reconstruct the scores of the individual matches.
The above-mentioned Maple package {\tt WorldCup} has a procedure, {\tt Khida}, that makes up such puzzles, and
another procedure, {\tt Sefer1}, that creates challenging puzzle books. For some samples, see

{\tt http://www.math.rutgers.edu/\~{}zeilberg/tokhniot/oWorldCupi} \quad, \quad
for $1 \leq i \leq 4$. Enjoy!

{\bf 7.} Happy birthday Richard, and keep up the good work! (and continue to practice what Phil Hanlon calls the ``the three H's'').

\bigskip
{\bf References}

[OEIS1] Neil J.A. Sloane, {\it Sequence A061927}, {\tt http://oeis.org/A061927}.

[EC1] Richard P. Stanley, ``{\it Enumerative Combinatorics, volume I}'',  first edition: Wadsworth \& Brooks/Cole, 1986;
second edition: Cambridge University Press, 2011. \hfill\break
Available on-line (viewed July 7, 2014): {\tt http://math.mit.edu/\~{}rstan/ec/ec1.pdf} \quad .
\bigskip
\hrule
\bigskip
Shalosh B. Ekhad, c/o D. Zeilberger, Department of Mathematics, Rutgers University (New Brunswick), Hill Center-Busch Campus, 110 Frelinghuysen
Rd., Piscataway, NJ 08854-8019, USA. 
\bigskip
Doron Zeilberger, Department of Mathematics, Rutgers University (New Brunswick), Hill Center-Busch Campus, 110 Frelinghuysen
Rd., Piscataway, NJ 08854-8019, USA. \hfill\break
url: {\tt http://www.math.rutgers.edu/\~{}zeilberg/}   
\quad . \hfill\break
Email: {\tt zeilberg at math dot rutgers dot edu}   \quad .
\bigskip
\hrule
\bigskip
{\bf EXCLUSIVELY PUBLISHED IN THE PERSONAL JOURNAL OF SHALOSH B. EKHAD and DORON ZEILBERGER} 
{\tt http://www.math.rutgers.edu/\~{}zeilberg/pj.html} and {\tt arxiv.org} \quad .
\bigskip
\hrule
\bigskip
July 7, 2014
\end